\documentclass[11pt]{amsart} 
\newtheorem{theorem}{Theorem}[section]
\newtheorem{corollary}[theorem]{Corollary}

\begin{document} 
\setlength{\unitlength}{0.01in}
\linethickness{0.01in}
\begin{center}
\begin{picture}(474,66)(0,0)
\multiput(0,66)(1,0){40}{\line(0,-1){24}}
\multiput(43,65)(1,-1){24}{\line(0,-1){40}}
\multiput(1,39)(1,-1){40}{\line(1,0){24}}
\multiput(70,2)(1,1){24}{\line(0,1){40}}
\multiput(72,0)(1,1){24}{\line(1,0){40}}
\multiput(97,66)(1,0){40}{\line(0,-1){40}}
\put(143,66){\makebox(0,0)[tl]{\footnotesize Proceedings of the Ninth Prague Topological Symposium}}
\put(143,50){\makebox(0,0)[tl]{\footnotesize Contributed papers from the symposium held in}}
\put(143,34){\makebox(0,0)[tl]{\footnotesize Prague, Czech Republic, August 19--25, 2001}}
\end{picture}
\end{center}
\vspace{0.25in}
\title[sublinear continuous order-preserving function]{Sublinear and 
continuous order-preserving functions for noncomplete preorders} 
\author{Gianni Bosi}
\address{Dipartimento di Matematica\\
Applicata {\em ``Bruno de Finetti''}\\ 
Universit\`a di Trieste, Piazzale\\
Europa 1, 34127, Trieste, Italy\\
Phone: +39 040 6767115. Fax: +39 040 54209.}
\email{giannibo@econ.univ.trieste.it}
\author{Magal\`i E. Zuanon}
\address{Istituto di Econometria e Matematica per le Decisioni Economiche\\
Universit\`a Cattolica del Sacro Cuore\\
Largo A. Gemelli 1, 20123, Milano, Italy}
\email{mzuanon@mi.unicatt.it}
\thanks{\textit{Corresponding author.} Gianni Bosi}
\thanks{Gianni Bosi and Magal\`i E. Zuanon, 
{\em Sublinear and continuous order-preserving functions for noncomplete 
preorders},
Proceedings of the Ninth Prague Topological Symposium (Prague, 2001),
pp.~1--8, Topology Atlas, Toronto, 2002}
\begin{abstract} 
We characterize the existence of a nonnegative, sublinear and continuous
order-preserving function for a not necessarily complete preorder on a
real convex cone in an arbitrary topological real vector space. 
As a corollary of the main result, we present necessary and sufficient
conditions for the existence of such an order-preserving function for a
complete preorder.
\end{abstract}
\maketitle

\section{Introduction} 

Necessary and sufficient conditions for the existence of a
continuous linear order-preserving function for a complete preorder on a
topological real vector space are already found in the literature (see
e.g.\ Candeal and Indur\'ain \cite{ci1995} and Neuefeind and Trockel
\cite{nt1995}). 
It is well known that there are important applications of such results in
expected utility theory and collective decision making. 
A characterization of the existence of a continuous linear utility
function for a complete preorder on a convex set in a normed real vector
space was presented by B\"{u}ltel \cite{b2001}. 
Further, some authors were concerned with the existence of a homogeneous
of degree one and continuous order-preserving function for a complete
preorder on a real cone in a topological real vector space (see e.g.\ Bosi
\cite{b1998}, Bosi, Candeal and Indur\'ain \cite{bci2000}, and Dow and
Werlang \cite{dw1992}). 
More recently, Bosi and Zuanon \cite{bz2001a} presented a characterization
of the existence of a nonnegative, homogeneous of degree one and
continuous order-preserving function for a noncomplete preorder on a real
cone in a topological real vector space. 
In a different context, other authors were concerned with the existence of
an additive order-preserving function on a completely preordered semigroup
(see e.g.\ Allevi and Zuanon \cite{az2000}, and Candeal, de Miguel and
Indur\'ain \cite{cdi1999}). 
Bosi and Zuanon \cite{bz2001b} presented a characterization of the
existence of a Choquet integral representation for a complete preorder on
the space of all continuous real-valued functions on a compact topological
space.

In this paper we provide an axiomatization of the existence
of a nonnegative, sublinear and continuous order-preserving function
for a not necessarily complete preorder on a real convex cone in a
topological real vector space. All the results we present are based
upon the notion of a {\em decreasing scale} (or {\em linear
separable system}) which was introduced by Herden \cite{h1989a, h1989b}
in order to characterize the existence of a continuous
order-preserving function for a preorder on a topological
space (see also Burgess and Fitzpatrick \cite{bf1977},
and Mehta \cite{m1998}). 

It should be noted that such a topic can be of some interest in the
applications to economics. 
Consider the following example concerning decision theory under
uncertainty.
Let ${\bf M}=\{\mu_{n}: n \in \{1,\ldots,n^{*}\}\}$ be
a finite family of {\em concave capacities} on
a measurable space $(\Omega , \mathcal{A})$, with $\Omega$
the {\em state space}, and $\mathcal{A}$ a {\em $\sigma$-algebra}
of subsets of $\Omega$. We recall that a {\em capacity} $\mu$
on $\mathcal{A}$ (i.e., a function from $\mathcal{A}$ into $[0,1]$ such
that 
$\mu ( \emptyset ) = 0$, 
$\mu ( \Omega ) = 1$, and 
$\mu (A) \leq \mu ( B)$ for all $A \subseteq B$, $A,B \in \mathcal{A}$) 
is said to be concave if for all sets
$A,B \in \mathcal{A}$, 
$$
\mu (A \cup B ) + \mu (A \cap B ) \leq \mu (A) + \mu (B)
$$
(see e.g.\ Chateauneuf \cite{c1996}). 
Let $X$ be a {\em real convex cone} of nonnegative {\em real random
variables} (i.e., measurable real functions) on $(\Omega , \mathcal{A})$,
and assume that $X$ is contained in $L^{1}(\Omega , \mathcal{A}, \mu_{n})$
for every $n \in \{1,\ldots,n^{*}\}$, where $L^{1}(\Omega, \mathcal{A}, \mu)$
stands for the normed space of all the real random variables $x$ such
that the {\em Choquet integral}
$$
\int_{\Omega}x d\mu = 
\int_{0}^{\infty}\mu(\{x \geq t\})dt + 
\int_{- \infty }^{0}(\mu(\{x \geq t\})-1)dt
$$
is finite (see e.g.\ Denneberg \cite{d1994}). 
Define a binary relation
$\preceq$ on $X$ as follows:
$$
x \preceq y \mbox{ if and only if }
\int_{\Omega}x d\mu_{n} \leq \int_{\Omega}y d\mu_{n}
\mbox{ for all }
n \in \{1,\ldots,n^{*}\}.
$$
It is clear that $\preceq$ is a preorder on $X$, and that $\preceq$ is not
complete in general. 
For every $n \in \{1,\ldots,n^{*}\}$, denote by $\tau_{n}$
the norm topology on $X$ which is associated to $\mu_{n}$, and let
$\tau$ be any (vector) topology on $X$ which is stronger than
$\tau_{n}$ for all $n \in \{1,\ldots,n^{*}\}$ (i.e., 
$\tau_{n} \subseteq \tau$ for all $n \in \{1,\ldots,n^{*}\}$). 
Then the real-valued function $u$ on $X$ defined by 
$$
u(x) =
\sum_{n=1}^{n^{*}}
\int_{\Omega}x d\mu_{n} \quad(x \in X)
$$
is a nonnegative, sublinear and $\tau$-continuous {\em order-preserving
function} for $\preceq$. 
Indeed, it is clear that $x \preceq y$ implies $u(x) \leq u(y)$ for all
$x,y \in X$.
If $x \prec y$, then we have $u(x) < u(y)$, since 
$\int_{\Omega}x d\mu_{n} \leq \int_{\Omega}y d\mu_{n}$ 
for all $n \in \{1,\ldots,n^{*}\}$,
and there exists at least one index $\bar n \in \{1,\ldots,n^{*}\}$ such
that $\int_{\Omega}x d\mu_{\bar n} < \int_{\Omega}y d\mu_{\bar n}$.
Further, $u$ is sublinear since the functional
$x \rightarrow \int_{\Omega}x d\mu_{n}$ is sublinear for all
$n \in \{1,\ldots,n^{*}\}$. 
Finally, $u$ is $\tau$-continuous since the functional 
$x \rightarrow \int_{\Omega}x d\mu_{n}$ is $\tau_{n}$-continuous, and
therefore $\tau$-continuous for all $n \in \{1,\ldots,n^{*}\}$ (see
Denneberg \cite[Proposition 9.4]{d1994}).

\section{Notation and preliminaries}

A {\em preorder}
$\preceq$ on an arbitrary set $X$ is a reflexive and transitive binary
relation on $X$. The {\em strict part} and the {\em symmetric part}
of a given preorder $\preceq$ will be denoted by $\prec$, and
respectively $\sim$. A preorder $\preceq$ on a set $X$ is said to be
{\em complete} if for any two elements $x,y \in X$ either $x \preceq y$
or $y \preceq x$. 

If $\preceq$ is a preorder on a set $X$,
then the pair $(X,\preceq )$ will be referred to as a {\em preordered set}.
Define, for every $x \in X$, $L_{\prec}(x)=\{z \in X: z \prec x\}$,
$U_{\prec}(x)=\{z \in X: x \prec z\}$. 

Given a preordered set
$(X, \preceq)$, a real-valued function $u$ on $X$ is
said to be
\begin{enumerate}
\item
{\em increasing} if $u(x) \leq u(y)$ for every $x,y \in X$ such that 
$x \preceq y$;
\item
{\em order-preserving} if it is increasing and $u(x) < u(y)$ for
every $x,y \in X$ such that $x \prec y$.
\end{enumerate}

If $(X, \preceq )$ is a preordered set, and $\tau$ is a topology on $X$,
then the triple $(X, \tau , \preceq )$ will be referred to as a {\em
topological preordered space}.
If $(X, \tau , \preceq )$ is a topological completely preordered space,
then the complete preorder $\preceq$ is said to be {\em continuous} if
$L_{\prec}(x)$ and $U_{\prec}(x)$ are open subsets of $X$ for every
$x \in X$.

Given a preordered set $(X, \preceq )$,
a subset $A$ of $X$ is said to be
{\em decreasing} if $y\in A$ whenever $y \preceq x$ and
$x \in A$.

In the sequel, the
symbol ${\mathbb Q}^{++}$ (${\mathbb R}^{++}$) will stand for the set of all
positive rational (real) numbers. If $(X, \tau )$ is a topological space,
then denote by $\overline A$ the topological closure of any
subset $A$ of $X$.

We say that a family
$\mathcal{G}=\{G_{r}: r \in {\mathbb Q}^{++}\}$ is a {\em countable
decreasing scale} ({\em countable linear separable system}) in a
topological preordered space $(X,\tau , \preceq )$ if
\begin{enumerate}
\item
$G_{r}$ is an open decreasing subset of $X$ for every $r \in
{\mathbb Q}^{++}$;
\item
$\overline{G_{r_{1}}} \subseteq G_{r_{2}}$ for every $r_{1},r_{2}
\in {\mathbb Q}^{++}$ such that $r_{1} < r_{2}$;
\item
$\displaystyle{\bigcup_{r \in {\mathbb Q}^{++}}G_{r}=X}$.
\end{enumerate}

If $E$ is a real vector space, then define, for every subset $A$ of $E$
and any real number $t$, $tA=\{ta:a \in A\}$. Further, if $A$ and $B$ are
any two subsets of a (real) vector space $E$, then define
$A+B=\{a+b:a \in A,b\in B\}$. 

A subset $X$ of a
{\em real vector space} $E$ is said to be
\begin{enumerate}
\item
a {\em real cone} if $tx \in X$ for every $x \in X$ and $t \in
{\mathbb R}^{++}$;
\item
a {\em real convex cone} if it is a real cone and $x+y \in X$ for every
$x,y \in X$.
\end{enumerate}

A real-valued function $u$ on a real cone $X$ in a real vector
space $E$ is said to be {\em homogeneous of degree one} if $u(tx) = t
u(x)$ for every $x \in X$ and $t \in {\mathbb R}^{++}$.

A real-valued function $u$ on a real convex cone $X$ in a real
vector space $E$ is said to be {\em sublinear} if it is homogeneous of
degree one and {\em subadditive} (i.e., $u(x+y) \leq u(x) + u(y)$ for
every $x,y \in X$).

Given a {\em topological real vector space} $E$, denote by $\tau$ the
vector topology for $E$ (i.e., the topology on $E$ which makes the vector
operations
continuous). 

If $X$ is any subset of a topological real vector space $E$, denote by
$\tau_{X}$ the topology induced on $X$ by the vector topology $\tau$ on
$E$.

If $(X, \preceq )$ is a preordered real cone in a topological
real vector space $E$, then we say that a countable
decreasing scale $\mathcal{G}=\{G_{r}: r \in {\mathbb Q}^{++}\}$
in $(X, \tau_{X},\preceq )$ is {\em homogeneous} if $qG_{r}=G_{qr}$ for
every $q,r \in {\mathbb Q}^{++}$.

If $(X, \preceq )$ is a preordered real convex cone in a topological real
vector space $E$, then we say that a countable decreasing scale
$\mathcal{G}=\{G_{r}: r \in {\mathbb Q}^{++}\}$ in $(X, \tau_{X},\preceq
)$ is
{\em subadditive} if $G_{q}+G_{r}\subseteq G_{q+r}$ for
every $q,r \in {\mathbb Q}^{++}$.

\section{Existence of a sublinear continuous order-preserving function}

In the following theorem we characterize the existence of a nonnegative,
sublinear and continuous order-preserving function for a not necessarily
complete preorder
on a real convex cone in a topological real vector space.

\begin{theorem}\label{thm}
Let $\preceq$ be a preorder on a real convex cone $X$ in a topological
real vector space $E$. 
Then the following conditions are equivalent:
\begin{enumerate}
\item\label{theorem1}
There exists a nonnegative, sublinear and continuous order-preserving
function $u$ for $\preceq$.
\item\label{theorem2}
There exists a homogeneous and subadditive countable decreasing scale
$\mathcal{G}=\{G_{r}: r \in {\mathbb Q}^{++}\}$ in $(X, \tau_{X},{\preceq})$ 
such that for every $x,y \in X$ with $x \prec y$ there exist 
$r_{1},r_{2} \in {\mathbb Q}^{++}$ with 
$$
r_{1} < r_{2}, 
x \in G_{r_{1}}\ \mbox{and}\
y \not \in G_{r_{2}}.
$$
\end{enumerate}
\end{theorem}

\begin{proof}
(\ref{theorem1}) $\Rightarrow$ (\ref{theorem2}).
Assume that there exists a nonnegative, sublinear and continuous
order-preserving function $u$ for $\preceq$. Define $G_{r}=u^{-1}([0,r[)$
for every $r \in {\mathbb Q}^{++}$. Let us show that
$\mathcal{G}=\{G_{r}: r \in {\mathbb Q}^{++}\}$ is a homogeneous and
subadditive countable decreasing scale satisfying condition
(\ref{theorem2}). 
Using the fact that $u$ is nonnegative and order-preserving, we have that 
for every $x,y \in X$ such that $x \prec y$, there exist
$r_{1},r_{2} \in {\mathbb Q}^{++}$ such that $u(x) < r_{1} < r_{2} < u(y)$, 
and therefore $x \in G_{r_{1}}$, $y \not \in G_{r_{2}}$. 
Further, since $u$ is homogeneous of degree one, 
$$
qG_{r}= qu^{-1}([0,r[) = u^{-1}([0,qr[) =G_{qr}\
\mbox{for every $q,r \in {\mathbb Q}^{++}$.}
$$ 
Hence, $\mathcal{G}$ is homogeneous. 
It only remains to show that $\mathcal{G}$ is subadditive. 
To this aim, consider any two rational numbers $q,r \in {\mathbb Q}^{++}$,
and let $z \in G_{q} + G_{r}$.
Then there exist two elements $x,y \in X$ such that
$z = x+y$, $u(x) < q$, $u(y) < r$. Hence, using the fact that $u$ is
subadditive, we have $u(z) = u(x+y) \leq u(x) + u(y) < q+r$, and
therefore $z \in G_{q+r}$.

(\ref{theorem2}) $\Rightarrow$ (\ref{theorem1}).
Define, for every
$x \in X$, \[u(x)= \inf \{r \in {\mathbb Q}^{++}: x \in G_{r}\}.\]
Then $u$ is a nonnegative continuous order-preserving function for
$\preceq$. 
Indeed, by using the fact that $G_{r}$ is a decreasing set for every 
$r \in {\mathbb Q}^{++}$, it is easily seen that $u$ is increasing. 
Further, $u$ is order-preserving by condition (\ref{theorem2}) above,
and $u$ is continuous since $u(x)= \inf \{r \in {\mathbb Q}^{++}: x \in
\overline{G_{r}}\}$ (see e.g.\ Theorem 1 in Bosi and Mehta
\cite{bm2001} for details). 
In order to show that $u$ is homogeneous of degree one, it suffices to
prove that for no $r \in {\mathbb Q}^{++}$, and $x \in X$ it is
$u(rx) \neq ru(x)$. 
Then the thesis follows by a standard continuity argument. 
This part of the proof is already found in Bosi and Zuanon 
\cite[Theorem 1]{bz2001a}.
Nevertheless, we present all the details here for reader's convenience.
By contradiction, assume that there exist $r \in {\mathbb Q}^{++}$,
and $x \in X$ such that $u(rx) < ru(x)$. 
Then, from the definition of $u$, there exists $r' \in {\mathbb Q}^{++}$
such that $u(rx) < r' < ru(x)$, $rx \in G_{r'}$. 
Since $u(x) > \frac{r'}{r}$, it follows that 
$x \not \in G_{\frac{r'}{r}}=\frac{1}{r}G_{r'}$, and therefore we arrive
at the contradiction $rx \not \in G_{r'}$. 
Analogously it can be shown that for no $r \in {\mathbb Q}^{++}$, and 
$x \in X$ it is $ru(x) < u(rx)$. 
It remains to prove that $u$ is subadditive.
By contradiction, assume that there exist two elements $x,y \in X$ such 
that $u(x) + u(y) < u(x+y)$. Then, from the definition of $u$, there exist 
two rational numbers $q,r \in {\mathbb Q}^{++}$ such that
$u(x) + u(y) < q + r < u(x+y)$, $x \in G_{q}$, $y \in G_{r}$. 
Since the countable decreasing scale 
$\mathcal{G}= \{G_{r}: r \in {\mathbb Q}^{++}\}$ is subadditive, we have
$x + y \in G_{q+r}$, and therefore $q + r < u(x+y)$ is contradictory from
the definition of $u$. 
This consideration completes the proof.
\end{proof}

If $\preceq$ is a {\em homothetic} complete preorder on a real cone $X$ in
a real vector space $E$ (i.e., $x \preceq y$ entails 
$tx \preceq ty$ for every $x,y \in X$, and $t \in {\mathbb R}^{++}$), then
consider the following subcones of $X$: 
$$X_{0} = \{x \in X: x \sim tx 
\mbox{ for some positive real number } t \neq 1\};$$ 
$$X_{+} = \{x \in X: x \prec tx \mbox{ for some real number } t > 1\};$$
$$X_{-} = \{x \in X: tx \prec x \mbox{ for some real number } t > 1\}.$$

In the following corollary, we present a characterization of the
existence of a nonnegative, sublinear and continuous order-preserving
function for a complete preorder on a real convex cone in a topological
real vector space. We recall that, given a preordered set $(X, \preceq)$,
a subset $A$ of $X$ is said to be an {\em order-dense subset of}
$(X, \preceq)$ if for every $x,y \in X$ such that $x \prec y$ there
exists $a \in A$ such that $x \prec a \prec y$.

\begin{corollary}\label{cor}
Let $\preceq$ be a complete preorder on a real convex cone $X$ in a
topological real vector space $E$, and assume that $X_{0}$ and $X_{+}$
are both nonempty, while $X_{-}$ is empty.
Then the following conditions are equivalent:
\begin{enumerate}
\item\label{cor1}
There exists a nonnegative, sublinear and continuous 
order-preserving function $u$ for $\preceq$.
\item\label{cor2}
The following conditions are satisfied: 
\begin{enumerate}
\item\label{cora}
$\preceq$ is homothetic;
\item\label{corb}
$\preceq$ is continuous;
\item\label{corc}
The set $\{qx_{+}: q \in {\mathbb Q}^{++}\}$ is an order-dense subset
of $(X_{+}, {\preceq})$ for every element $x_{+} \in X_{+}$;
\item\label{cord}
$x \sim y$ for every $x,y \in X_{0}$;
\item\label{core}
$x \prec x_{+}$ for every $x \in X_{0}$, $x_{+} \in X_{+}$;
\item\label{corf}
$x + y \prec (q+r)x_{+}$ for every $x,y \in X$, $x_{+} \in X_{+}$, 
$q,r \in {\mathbb Q}^{++}$ such that $x \prec qx_{+}$,
$y \prec rx_{+}$.
\end{enumerate}
\end{enumerate}
\end{corollary}

\begin{proof}
(\ref{cor1}) $\Rightarrow$ (\ref{cor2}).
Assume that there exists a nonnegative, sublinear and continuous
order-preserving function $u$ for $\preceq$. 
Then it is clear that $\preceq$ is homothetic and continuous. 
If we consider any element $x_{+} \in X$ such that $x_{+} \prec tx_{+}$
for some real number $t>1$, then it is necessarily $u(x_{+}) > 0$, and
therefore, using the fact that $u$ is homogeneous of degree one, condition
(\ref{corc}) easily follows. 
Finally, it is easily seen that conditions (\ref{cord}), (\ref{core}) and
(\ref{corf}) are verified.

(\ref{cor2}) $\Rightarrow$ (\ref{cor1}). 
Consider any element $x_{+} \in X$ such that $x_{+} \prec tx_{+}$ for some
real number $t>1$, and let $G_{r}=L_{\prec}(rx_{+})$ for every 
$r \in {\mathbb Q}^{++}$. 
Then it is easy to check that the family 
$\mathcal{G}= \{G_{r}: r \in {\mathbb Q}^{++}\}$
is a countable decreasing scale satisfying condition
(\ref{theorem2}) of Theorem \ref{thm}. 
Indeed, by homotheticity of $\preceq$, $x \prec rx_{+}$ is
equivalent to $qx \prec qrx_{+}$ ($q \in {\mathbb Q}^{++}$), and therefore
$\mathcal{G}$ is homogeneous (see Bosi and Zuanon 
\cite[Corollary 2]{bz2001a}).
Further, $\mathcal{G}$ is subadditive by condition (\ref{corf}) since,
for every $q,r \in {\mathbb Q}^{++}$, and $x,y \in X$, $x \in G_{q}$ and
$y \in G_{r}$ is equivalent to $x \prec qx_{+}$ and $y \prec rx_{+}$, which
implies $x+y \prec (q+r)x_{+}$ or equivalently $x+y \in G_{q+r}$. 
Finally, by condition (\ref{corc}) above, for every $x,y \in X$ such that
$x \prec y$ there exist $r_{1},r_{2} \in {\mathbb Q}^{++}$ such that 
$$
r_{1} < r_{2}, \quad
x \prec r_{1}x_{+} \prec r_{2}x_{+} \prec y,
$$ 
or equivalently 
$$
x \in L_{\prec}(r_{1}x_{+}), \quad
y \not \in L_{\prec}(r_{2}x_{+}).
$$ 
So the proof is complete.
\end{proof}

Denote by $\bar 0$ the {\em zero vector} in a real vector space $E$. 
In the following corollary we are concerned with a sublinear 
representation of a complete preorder on a real convex cone containing the
zero vector.

\begin{corollary}
Let $\preceq$ be a complete preorder on a real convex cone $X$ in a 
topological real vector space $E$, and assume that $X_{+}$ is nonempty,
while $X_{-}$ is empty. 
If in addition $\bar 0$ belongs to $X$, then there exists a nonnegative,
sublinear and continuous order-preserving function $u$ for $\preceq$ if
and only if $\preceq$ is homothetic and continuous, and it satisfies
condition (\ref{corf}) of Corollary \ref{cor}.
\end{corollary}

\begin{proof}
From the corollary in Bosi, Candeal and Indur\'ain \cite{bci2000}, there
exists a nonnegative, homogeneous of degree one and continuous
order-preserving function $u$ for $\preceq$. 
Indeed, the complete preorder $\preceq$ on the real (convex) cone $X$ is
homothetic and continuous, and we have in addition $\bar 0 \in X$. 
Let us show that $u$ must be subadditive as a consequence of condition
(\ref{corf}) of Corollary \ref{cor}.
Assume by contraposition that there exist $x,y \in X$ such that 
$u(x) + u(y) < u(x+y)$, and consider any element $x_{+} \in X_{+}$.
Then it must be $u(x_{+}) > 0$ since $u$ is a homogeneous of degree one
utility function for $\preceq$, and there exist two positive rational
numbers $q$ and $r$ such that 
$$
u(x) + u(y) < (q+r)u(x_{+}) < u(x+y),\quad
x \prec qx_{+}, \quad
y \prec rx_{+}.
$$ 
But here we have a contradiction since it should be 
$u(x+y) < (q+r)u(x_{+})$ by condition (\ref{corf}) of Corollary \ref{cor}.
So the proof is complete.
\end{proof}

\providecommand{\bysame}{\leavevmode\hbox to3em{\hrulefill}\thinspace}
\providecommand{\MR}{\relax\ifhmode\unskip\space\fi MR }
\providecommand{\MRhref}[2]{%
  \href{http://www.ams.org/mathscinet-getitem?mr=#1}{#2}
}
\providecommand{\href}[2]{#2}

\end{document}